\documentclass{amsart}
\usepackage{amsmath, amsfonts, amsthm}
\usepackage[foot]{amsaddr}
\def\Iscr{\mathcal{I}}
\DeclareMathOperator{\cl}{cl}
\DeclareMathOperator{\rank}{rank}

\theoremstyle{plain}

  \theoremstyle{definition}
  
  \theoremstyle{definition}
  \newtheorem*{example*}{\protect\examplename}
  \theoremstyle{plain}
  
  \theoremstyle{plain}
  
  \theoremstyle{plain}
  
 \theoremstyle{definition}
 \newtheorem*{defn*}{\protect\definitionname}
  \theoremstyle{plain}
  \newtheorem*{lem*}{\protect\lemmaname}
  \theoremstyle{plain}
  \newtheorem*{thm*}{\protect\theoremname}  

\newtheorem{defin}{Definition}

\newtheorem{theo}{Theorem}
\newtheorem{rem}{Remark}
\mathchardef\cdotchar=12801

\title{A family of rank $4$ non-algebraic matroids with pseudomodular dual}
\author{Winfried Hochst\"attler}
\address{Winfried Hochst\"attler\\FernUniversit\"at in Hagen\\D-58084 Hagen\\Germany} 
\begin{document}
\begin{abstract}
  The Tic-Tac-Toe matroid is a paving matroid of rank $5$ on 9
  elements which is pseudomodular and whose dual is non-algebraic. It
  has been proposed as a possible example of an algebraic matroid
  whose dual is not algebraic.  We present an infinite family of
  matroids sharing these properties and generalizing the Tic-Tac-Toe matroid.
  
\end{abstract}
\maketitle
\section{Introduction}
Hassler Whitney~\cite{whitney35} introduced matroids 1935 as
structures that capture the ``abstract properties of linear
dependence''. In the second edition of his ``Moderne Algebra'' van
der Waerden~\cite{vanderWaerden} remarked 1937, that algebraic dependencies in field
extensions share these abstract properties, a fact that had already been
shown by Steinitz in 1910~\cite{Steinitz}. 

\begin{defin}
  Let $k$ be a field and $k\subseteq K$ a field extension. A finite
  subset $\{v_1,\ldots,v_k\} \subseteq K$ is called {\em algebraically
  dependent} over $k$ if there exists a non-trivial polynomial $0 \ne
  p(x_1,\ldots,x_k) \in k[x_1,\ldots,x_k]$ such that 
  $$p(v_1,\ldots,v_k)=0.$$
  The set $K$ is called {\em algebraically
    independent}, otherwise.
\end{defin}

Algebraic independence satisfies the axioms of matroid theory (see e.g.\ 
6.7.1 in~\cite{Oxley}).

\begin{defin}
  A matroid $M=(E,\mathcal{I})$ on a finite set $E$ is called {\em
    algebraic}, if there exists a field extension $k \subseteq K$ a
  subset $E' \subseteq K$ and a bijection $\sigma:E \to E'$ such that
  for all $I\subseteq E$
  $$I \in \mathcal{I} \Leftrightarrow \sigma(I)\mbox{ is algebraically independent over }k.$$
\end{defin} 
  
It has been shown by Ingleton\cite{Ingleton} only in 1971 that
algebraic matroids are a strict superclass of the class of linear
matroids and it took until 1975 for Ingleton and Main~\cite{InglMain}
to show that ``Non-algebraic matroids exist'', by proving that the
V\'amos matroid is non-algebraic.  For that purpose they showed that
in an algebraic matroid the three bounding lines of a prism must
intersect in a common point of the closure of the representing field.
We will give an exact proof in the next section.

Dress and Lov\'asz~\cite{DressLovasz} generalized this to the ``weak
series reduction theorem'' and Bj\"orner and Lov\'asz~\cite{BjoLo}
used this to define pseudomodular lattices.

The key lemma in the proof that the V\'amos matroid is non-algebraic
was generalized by A.~Dress and L.~Lov\'asz to the ``series reduction
theorem''~\cite{DressLovasz}.

Considering the question whether pseudomodularity of the geometric
lattice of a matroid is sufficient for the existence of a lattice
theoretical dual, an adjoint of that matroid, Alfter and
Hochst\"attler~\cite{alfter} presented a pseudomodular matroid of rank
$5$ which does not admit an adjoint and named it the Tic-Tac-Toe
matroid. It will appear in Section~\ref{sec:dual} as $M_3^\ast$. It had
been observed independently by M.~Alfter and by B.~Lindstr\"om that
the dual of this matroid is non-algebraic, using an application of the
key lemma of Ingleton and Main.

Algebraic matroids are closed under taking minors and under
truncation. It is still an open problem though, whether they are
closed under matroid duality.

Since pseudomodularity in a certain sense seems to capture the most
important properties of algebraic matroids, it was suggested by Lazlo
Lov\'asz \cite{dimacsfloor} that the Tic-Tac-Toe matroid might be a
good candidate as a counterexample for the  unsolved question
whether the dual of an algebraic matroid is algebraic again.

We tried to popularize this suggestion with a preprint~\cite{zpr}
in 1997 and a note in the Proceedings of the Graph Theory
V~\cite{gtv}. We were partially successfull. The author is aware of at
least two PhD projects which were motivated by the question, whether
the Tic-Tac-Toe matroid is algebraic or not, namely by Stephan
Kromberg~\cite{Kromberg} and Guus Bollen~\cite{Bollen}. Although both
theses present very nice results, they did not settle the original
question.

Algebraic matroids are closed under taking minors and under
truncation. It is still an open problem though, whether they are
closed under matroid duality.

The paper is organized as follows. In the next section we present the
combinatorial properties of algebraic matroids that led to the
definition of pseudomodularity. In Section~\ref{sec:primal} we define
an infinite family of rank $4$ paving matroids generalizing the dual
of the Tic-Tac-Toe matroid and prove that all of them are
non-algebraic. In Section~\ref{sec:dual} we show that all their duals
are pseudomodular, and we end with some remarks and open questions in
Section~\ref{sec:5}.  We assume some familiarity with matroid theory,
the standard reference is \cite{Oxley}.

\section{Combinatorial Properties of Algebraic Matroids}

\begin{defin}
  Let $M=(E,\Iscr)$ be a matroid on a finite set $E$ and $S \subseteq
  A \subseteq E$. Then $S$ is {\em in series in $A$} if contracting $A
  \setminus S$ turns $S$ into a circuit.
\end{defin}

\begin{theo}[Dress-Lov\'asz 1987]\label{DressLovasz}
  Let $M=(E,\Iscr)$ be an algebraic matroid represented by a set $E'
  \subseteq K$ over a field $k$, $S'$ in series in $A'\subseteq E'$ and $K$
  algebraically closed. Then there exists  $\beta \in K$ such that $\forall T' \subseteq A' \setminus S'$: 
  $$S' \cup T' \mbox{ is
    algebraically dependent} \Leftrightarrow \beta \cup T' \mbox{ is
    algebraically dependent}.$$
\end{theo}

Thus, in an algebraic matroid every series admits a shortcut.
From this it is seen that the V\'amos matroid is non-algebraic as
follows: 

The V\'amos matroid $V$ is a paving matroid, i.e.\ all its circuit
have $\rank(V)$ or $\rank(V)+1$ elements,  on an eight point
set $\{a,a',b,b',c,c',d,d'\}$. Its five $4$-element circuits or {\em circuit hyperplanes} are given by
$\left\{a,a',b,b'\right\}$, $
  \left\{a,a',c,c'\right\}$, $\left\{b,b',c,c'\right\}$, $ \left\{b,b',d,d'\right\}$ and $\left\{c,c',d,d'\right\}.$
In particular  $\{a,a',d,d'\}$ is independent.  Assume $V$
were algebraically represented over a field extension $k \subseteq K$.
Since $\{a,a'\}$ is in series in $\{a,a',b,b',c,c'\}$ there exists a
$\beta_1$ in the algebraic closure of $K$ which lies on the
intersection of the two ``lines'' $\beta_1\in \overline{b,b'}
\,\,\cap\,\, \overline{c,c'}$. We denote the closure of $V \cup
\beta_1$ by overlining the sets. We also have $\beta_1 \in
\overline{a,a',b,b'} \,\,\cap\,\, \overline{a,a',c,c'}=
\overline{a,a'}$. Hence, according to the series reduction theorem the three bounding lines of the ``prism'' ${a,a',b,b',c,c'}$ must intersect in a point in the algebraic closure of the representing field.

By symmetry, there also exists $\beta_2\in
\overline{a,a'}\,\,\cap\,\,\overline{b,b'} \,\,\cap\,\,
\overline{d,d'}$. Since $\beta_1$ and $\beta_2$ lie in the
intersection of the same lines they must be parallel. Thus the two
lines $\overline{a,a'}$ and $\overline{d,d'}$ have an intersection of
rank at least one, contradicting the independence of the set
$\{a,a',d,d'\}$.
    
    All proofs of non-algebraicity, known to the author, apply the
    series reduction theorem, postulate additional points and derive a
    contradiction. 

\section{The Primal Matroids}\label{sec:primal}
For ease of definition we start with the generalizations of the dual
of the Tic-Tac-Toe Matroid which all are non algebraic.
\begin{defin}
  Let $k\ge 3$ and $E=\{a_1,\ldots, a_k, b_1,\ldots,b_k, c_1, \ldots,
  c_k\}$. We define a paving matroid (see~\cite{Oxley}, Proposition
  1.3.10) $M_k$ of rank 4 on $E$ by defining a set of circuit
  hyperplanes. These are given by
\[
\{(a_i,a_{j},b_i,b_{j})\mid i < j\} \cup \{(a_i,a_{j},c_i,c_{j})\mid i < j, (i,j) \ne (1,k)  \} \cup  
\{(b_i,b_{j},c_i,c_{j})\mid i \ne j\}.
\]

Note, that $\{(a_1,a_k,c_1,c_k)\}$ is a basis and not a circuit hyperplane.
\end{defin}

\begin{theo}
  $M_k$ is not algebraic.
\end{theo}
\begin{proof}
  $M_k$ may be considered as a complete graph, where each edge has
  been replaced by a prism, but one hyperplane $\{(a_1,a_k,c_1,c_k)\}$
  is broken. In particular the upper edges of the ``triangle''
  $(k-1,k,1)$, namely $\{(a_{k-1},a_{k},b_{k-1},b_{k})\}$, $
\{(a_{k-1},a_{k},c_{k-1 },c_{k})\}$ and
$\{(b_{k-1},b_k,c_{k-1},c_k)\}$ as well as
$\{(a_{1},a_{k},b_{1},b_{k})\}$, \linebreak[4] $\{(a_{1},a_{{k}},c_{1},c_{{k}})\}$ and
 $\{(b_{1},b_{k},c_{1},c_{k})\}$ form such a prism. Note
that the ``prism'' of the edge $(1,k)$ is the one with the broken
hyperplane.

  Assume $M_k$ were algebraic.  According to the Ingleton-Main Lemma (compare the proof of non-algebraicity of the V\'amos-matroid) 
  there exist
  \begin{eqnarray*}
&p \in \cl(a_{k-1},a_k) \cap \cl(b_{k-1},b_k) \cap
  \cl(c_{k-1},c_k) \text{ and } \\ &p' \in \cl(a_{1},a_{k-1}) \cap
  \cl(b_{1},b_{k-1}) \cap \cl(c_{1},c_{k-1}).     
  \end{eqnarray*}
 Hence $\{p,p'\}$ is a
  subset of $\cl(a_{k-1},a_k,a_1)$ as well as of
  $\cl(b_{k-1},b_k,b_1)$ and $\cl(c_{k-1},c_k,c_1)$. This implies that the eight
  points $\{a_1,b_1,c_1,a_k,b_k,c_k,p,p'\}$ form a Vamos matroid,
  contradicting algebraic representability.
\end{proof}

\begin{rem}
  The choice of the third vertex $k-1$ in the triangle is arbitrary. It
  could be replaced by any vertex from $2$ to $k-1$ to
  give rise to a copy of $M_3$, the dual of the Tic-Tac-Toe matroid.
\end{rem}

\section{The Dual Matroids}\label{sec:dual}
Since $M_k$ is sparse paving, so is its dual $M_k^\ast$. Note, that
its circuit hyperplanes are given by the complements of the circuit
hyperplanes of $M_k$.

We will show that $M_k^\ast$ is pseudomodular for all $k \ge 3$. Since
we have to deal with both matroid theoretic as well as lattice
theoretical concepts we will frequently identify elements of the
geometric lattice with the set of elements of the matroid in the
corresponding flat.

\begin{defin}\label{def:pm}
  Let $M$ be a matroid and $L$ be its geometric lattice of closed
  flats. Then $M$ is {\em pseudomodular}~\cite{BjoLo}, if $\forall a,b,c \in L:$
  \begin{eqnarray}
    \label{eq:1}
&&  r(a\vee b \vee c) - r(a\vee b)=r(a \vee c) - r(a)=r(b \vee c) - r(b)\\&\Longrightarrow&
r((a\vee c) \wedge (b \vee c)) - r(a\wedge b)=r(a \vee c) - r(a).      \label{eq:1a}
  \end{eqnarray}
  
\end{defin}
\begin{theo}
  $M_k^\ast$ is pseudomodular.
\end{theo}

\begin{proof}
  Assume $M_k^\ast$ were not pseudomodular and let $a,b,c \in L$
  violate Definition~\ref{def:pm}. 
Since $r(a\vee c)+r(a\vee
  b)=r(a\vee b \vee c) + r(a)$ the flats $a\vee c$ and $a\vee b$ must form a
  modular pair with meet $a$. By symmetry we also have $(b \vee c)
  \wedge (a \vee b)= b$ and hence 
  \begin{equation}
    \label{eq:3}
  (a\vee c) \wedge (b \vee c) \wedge (a \vee b)= a \wedge b.
  \end{equation}
  By submodularity we have 
  \begin{eqnarray*}
  &&  r((a \vee c) \wedge (b \vee c))+r( a \vee b) \\&\ge & r((a \vee c) \wedge (b \vee c))\vee (a \vee b)) + r((a \vee c) \wedge (b \vee c)\wedge (a \vee b)) \\&=& r(a \vee b \vee c) + r(a\wedge b).
  \end{eqnarray*}
Thus, violating Definition~\ref{def:pm}  implies 
    \begin{equation}
      \label{eq:2}
    r((a\vee c) \wedge (b \vee c)) - r(a\wedge b)\ge r(a\vee b \vee c)-r(a\vee b)+1=r(a \vee c) - r(a) +1.        
    \end{equation}
  By  definition of $M_k^\ast$ any flat $F$ with $r(F) \le r(M_k^\ast)-2$
  is independent, implying $r(F)=|F|$.
Now
\begin{eqnarray*}
  |a\vee c|&=&|(a\vee c)\wedge (a \vee b)| + |(a\vee c)\wedge (b \vee c)|-|(a\vee c) \wedge (b \vee c) \wedge (a \vee b)|\\
  & \stackrel{\eqref{eq:3}}{=}&r((a\vee c)\wedge (a \vee b)) + r((a\vee c)\wedge (b \vee c))-r(a \wedge b)\\
  &\stackrel{\eqref{eq:2}}{\ge}&r((a\vee c)\wedge (a \vee b)) + r(a\vee c) -r(a)+1
  \\&\ge&r(a) + r(a\vee c) -r(a)+1\\
  &=& r(a\vee c)+1. \end{eqnarray*}
Hence, $(a\vee c)$ must be a circuit hyperplane of $M_k^\ast$ and by symmetry this also holds for $(b\vee c)$. 
A similar computation shows that the same is true for $a \vee b$:
\begin{eqnarray*}
  |a\vee b|&=&r((a\vee b)\wedge (b \vee c)) + r((a\vee b)\wedge (a \vee c))-r(a \wedge b)\\
&=& r(b) + r(a) - r(a \wedge b)\\
&\stackrel{\eqref{eq:2}}{\ge}&r(b) + r(a \vee c) - r((a\vee c)\wedge (b \vee c)) + 1\\
&\stackrel{\text{subm.}}{\ge}&r(b) + r(a \vee c) - (r((a\vee c)+ r(b \vee c))- r(a\vee b \vee c))+1\\
&=&r(b)-r(b \vee c) + r( a \vee b \vee c)+1\\
&\stackrel{\eqref{eq:1}}{=}& r(a\vee b)+1. \end{eqnarray*}

Using the defining properties of $a,b,c$ this moreover implies
\[r(a)=r(b)=r(a \vee b \vee c) - 2.\]

Now, if $r((a\vee c) \wedge (b \vee c)) \le r(a\vee b \vee c) -3$ we conclude
 $r(a \wedge b)\le r(a\vee b \vee c) -5$ and thus
\begin{eqnarray*}
  |a\vee b| &=& r(a) + r(b) -r(a \wedge b)\\
 &\stackrel{\eqref{eq:2}}{\ge}&r(a) + r(b)) - r(a\vee b\vee c)+5\\           
 &=&r(a\vee b\vee c))+1,           
\end{eqnarray*} contradicting $|a \vee b|=r(a\vee b\vee c)$ as $a \vee b$ is a circuit hyperplane.

Hence, $a\vee b, a \vee c$ and $b\vee c$ are three circuit hyperplanes
which pairwise intersect in the three different colines $a,b$ and
$(a\vee c) \wedge (b \vee c)$.  But say $E \setminus \{a_i,a_j, b_i,
b_j\}$ intersects only two other circuit hyperplanes, namely $E
\setminus \{a_i,a_j, c_i, c_j\}$ and $E\setminus \{b_i,b_j, c_i,
c_j\}$ (provided $(i,j)\ne (1,k))$, in the coline $E\setminus
\{a_i,a_j, b_i,b_j, c_i, c_j\}$. Hence, the three circuit hyperplanes
could pairwise  intersect only in the same coline, a contradiction.

We conclude that $M_k^\ast$ must be pseudomodular.
\end{proof}

\section{Remarks and Open Problems}\label{sec:5}
Stephan Kromberg~\cite{Kromberg} analyzes two matroids which differ
only little from the Tic-Tac-Toe matroid, one of them linear, the
other one non-algebraic. 

Bamiloshin et al.\ \cite{bamiloshin} recently
enumerated all rank $4$ matroids on $8$ elements that are not linearly
representable and found $3$ examples, where they could not decide
whether they are algebraic or not. They furthermore described several
paving matroids of rank $5$ which are related to the Tic-Tac-Toe
Matroid and refer to~\cite{Bollen} concerning their representability
status.

It is not difficult to construct a pseudomodular matroid which contracts to a
non-algebraic matroid, and thus is non-algebraic itself. It might be
worthwhile to examine whether  some of our matroids
$M_k^\ast$ are non-algebraic for some $k\ge 4$.

The non-algebraicity of the matroids $M_k$ is not very exciting, as all
of them contain $M_3$. We wonder whether our construction of replacing
edges of a graph by prisms could avoid that, i.e.\ does there exist
a triangle-free two-connected graph, such that replacing its edges by
prisms and relaxing (breaking) one circuit hyperplane yields a
non-algebraic matroid?

 \bibliographystyle{amsplain} \bibliography{tictacplus}

\end{document}